\documentclass[12pt]{amsart}

\newtheorem{thm}{Theorem}[section]
\numberwithin{equation}{section}

\newtheorem{prop}[thm]{Proposition}
\newtheorem{cor}[thm]{Corollary}

\newtheorem{remark}[thm]{Remark}

\newcommand{\K}{\mathcal{K}}

\newcommand{\set}[1]{\left\{#1\right\}}

\newcommand{\norm}[1]{\left\Vert#1\right\Vert}

 \newcommand{\hot}{\hat{\otimes}}
 \newcommand{\al}{\alpha}
 \newcommand{\ep}{\varepsilon}

\begin{document}

\title[Banach algebra valued functions]{Amenability properties of Banach algebra valued continuous functions}

\author[R. Ghamarshoushtari]{Reza Ghamarshoushtari\dag}

\address{\dag\; Department of Mathematics\\ University of Manitoba\\ Winnipeg, Manitoba\\ R3T~2N2, Canada}

\email{\dag \; umghamar@cc.umanitoba.ca}

\author[Y. Zhang]{ Yong Zhang \ddag}
\address{\ddag \; Department of Mathematics\\
           University of Manitoba\\
           Winnipeg, Manitoba\\
           R3T 2N2 Canada}
\email{\ddag \; zhangy@cc.umanitoba.ca}
\thanks{\ddag \; Supported by NSERC Grant 238949-2005}




\begin{abstract}
  Let $X$ be a compact Hausdorff space and $A$ a Banach algebra. We investigate amenability properties of the algebra $C(X,A)$ of all $A$-valued continuous functions. We show that $C(X,A)$ has a bounded approximate diagonal if and only if $A$ has a bounded approximate diagonal; if $A$ has a compactly central approximate diagonal (unbounded) then $C(X,A)$ has a compactly approximate diagonal. Weak amenability of $C(X,A)$ for commutative $A$ is also considered.
\end{abstract}

\maketitle

\section{Introduction}

Let $A$ be a (complex) Banach algebra and $X$ a compact Hausdorff space. The Space of $A$-valued continuous functions on $X$ is denoted by $C(X,A)$. With pointwise algebraic operations and the uniform norm 
\[
\|f\|_\infty = \sup \{\|f(x)\|_A :\, x\in X\} \quad (f\in C(X,A)),
\]
$C(X, A)$ is a Banach algebra, where $\|\cdot \|_A$ denotes the norm of $A$. The study of the algebra $C(X, A)$ was initiated by I. Kaplansky in \cite{Kapl}, where he studied the structure of ideals of $C(X, A)$, in particular for $X$ being totally disconnected. Later  A. Hausner showed in \cite{Hausn} that the maximal ideal space of the general $C(X,A)$ is homeomorphic with $X\times \mathfrak{M}(A)$, where $\mathfrak{M}(A)$ is the maximal ideal space of $A$.

We study amenability properties of the Banach algebra $C(X,A)$. The notion of amenability for Banach algebras was first introduced by B. E. Johnson in 1972 in ~\cite{Johnson1}, where he showed that the group algebra $L^1(G)$ of a locally compact group $G$ is amenable as a Banach algebra if and only if the corresponding group $G$ is an amenable group. later, using Johnson's result on $L^1(G)$ and the Stone-Weierstrass Theorem M. V. \v{S}e\v{\i}nberg showed that $C(X)$ ($= C(X, \mathbb{C})$, the algebra of complex-valued continuous functions) is amenable for any compact Hausdorff space $X$ \cite{Sheinberg}. A constructive proof for this remarkable result was given by M. Abtahi and the second author in \cite{A-Z}. Here we note that, unlike $C(X)$, $C(X,A)$ is in general not a $C^*$-algebra and is no longer commutative if $A$ is not so.
However the method of \cite{A-Z} can be modified to deal with general $C(X,A)$. We shall show in Section~\ref{amenable} that $C(X,A)$ is amenable if and only if the range algebra $A$ is amenable. The proof uses significantly the following Grothendieck's inequality \cite[Theorem 5.5]{Pisier}.

\begin{thm}[Grothendieck]\label{Grothend}
Let $K_1, K_2$ be compact Hausdorff spaces, and let $\Phi$ be a bounded scalar-valued bilinear form on $C(K_1)\times C(K_2)$. Then there are probability measures $\mu_1, \mu_2$ on $K_1, K_2$, respectively, and a constant $k>0$ such that 
\[
 |\Phi(x_1,x_2)| \leq k \|\Phi\|\left(\int_{K_1}{|x_1|^2}d\mu_1 \int_{K_2}{|x_2|^2}d\mu_2\right)^\frac{1}{2} \quad(x_1\in C(K_1), x_2\in C(K_2)).
\]
\end{thm}

The smallest constant $k$ in the above theorem is called the Grothendieck constant, denoted $K_G^\mathbb{C}$. We know that $4/\pi \leq K_G^\mathbb{C} < 1.405$ \cite{Haag}. 

Let $A,B$ be two Banach spaces. We denote the Banach space projective tensor product of $A$ and $B$ by $A\hot B$, which is the completion of the algebraic tensor product $A\otimes B$ with respect to the projective norm
\[
\|u\|_p = \inf \set{\sum_{i=1}^{n}{\|a_i\| \|b_i\|}:\, u= \sum_{i=1}^{n}{a_i\otimes b_i} \in A\otimes B},
\]
where the infimum is taken over all representatives of $u$. The dual space of $A\hot B$ is identical to $BL(A,B;\mathbb{C})$, the space of all bounded scalar-valued bilinear forms on $A\times B$.

The Grothendieck Theorem stated above yields the following important inequality.

\begin{cor}\label{ineq}
Let $K_1, K_2$ be compact Hausdorff spaces. Then for each $u = \sum_{i=1}^{n}{a_i\otimes b_i} \in C(K_1)\otimes C(K_2)$ we have 
\[
\|u\|_p \leq c \left(\norm{\sum_{i=1}^{n}{|a_i|^2}}_\infty + \norm{\sum_{i=1}^{n}{|b_i|^2}}_\infty\right)
\]
where $c = \frac{1}{2}K_G^\mathbb{C}$. 
\end{cor}

\begin{proof}
From Theorem~\ref{Grothend} we have
\begin{align*}
\|a\otimes b\|_p & = \sup_{\Phi\in BL(A,B;\,\mathbb C)}|\Phi(a,b)| \\
                 &\leq K_G^\mathbb{C} \left(\int{|a|^2}d\mu_1 \int{|b|^2}d\mu_2\right)^\frac{1}{2}\leq c\left(\int{|a|^2}d\mu_1 + \int{|b|^2}d\mu_2\right)
\end{align*}
for $a\in C(K_1), b\in C(K_2)$. Thus
\[
\|u\|_p \leq c\left(\int{\sum_{i=1}^{n}{|a_i|^2}} d\mu_1 + \int{\sum_{i=1}^{n}{|b_i|^2}}d\mu_2\right)\leq c \left(\norm{\sum_{i=1}^{n}{|a_i|^2}}_\infty + \norm{\sum_{i=1}^{n}{|b_i|^2}}_\infty\right)
\]
for $u = \sum_{i=1}^{n}{a_i\otimes b_i} \in C(K_1)\otimes C(K_2)$.

\end{proof}

Let $A$ be a Banach algebra. Then $A\hot A$ is naturally a Banach $A$-bimodule. A net $(u_\al)\subset A\hot A$ is an \emph{approximate diagonal} for $A$ if
\[
a\cdot u_\al - u_\al\cdot a \to 0, \quad \pi(u_\al)a \to a \quad \text{for all } a\in A,
\]
where $\pi$: $A\hot A \to A$ is the product map defined by $\pi(a\otimes b) = ab$. The approximate diagonal $(u_\al)$ is called \emph{bounded} if it is a norm bounded net. It is called \emph{central} if $a\cdot u_\al = u_\al\cdot a$ for all $a\in A$ and all $\al$. The approximate diagonal $(u_\al)$ is a \emph{compactly approximate diagonal} for $A$ if for each compact set $K\subset A$ and any $\ep>0$ there is $\nu_0$ such that $\|a\cdot \al_\nu - \al_\nu\cdot a\|_p < \ep$ and $\|\pi(\al_\nu)a -a\| < \ep$ for all $a\in K$ whenever $\nu \geq \nu_0$. 

 B. E. Johnson showed in \cite{Johnson2} that a Banach algebra $A$ is amenable if and only if there is a bounded approximate diagonal for $A$. Our proof for the amenability of $C(X,A)$ in Section~\ref{amenable} is essentially based on this characterization of amenability -- we construct a bounded approximate diagonal for $C(X,A)$. Unbounded approximate diagonals for various Banach algebras have been studied in \cite{G-Z, C-G-Z}. Using the argument of \cite[Proposition 4.4]{G-Z} one can see that every Segal algebra on a [SIN] group has a compactly approximate diagonal. It is also easy to see that a $c_0$ or an $\ell_p$($p\geq 1$) direct sum of contractible Banach algebras has a central compactly approximate diagonal. We will show in Section~\ref{unbounded a.d.} that $C(X,A)$ has a compactly approximate diagonal if $A$ has a central compactly approximate diagonal. We will also show that $C(X,A)$ is weakly amenable if $A$ is a commutative weakly amenable Banach algebra with a bounded approximate identity.

\section{Amenable $C(X,A)$}\label{amenable}

\begin{thm}\label{bounded}
Let $X$ be a compact Hausdorff space and let $A$ be a Banach algebra. If $A$ has a bounded approximate diagonal, then so does $C(X,A)$.
\end{thm}

\begin{proof}
For $u = \sum_i{u_i\otimes v_i}\in C(X)\otimes C(X)$ and $\al = \sum_j{\al_j\otimes \beta_j} \in A\hot A$, it is readily seen that
\[
T(u,\al) = \sum_{i,j}u_i\al_j \otimes v_i\beta_j \in C(X,A)\hot C(X,A)
\]
and
\begin{equation}\label{mixed tensor}
\|T(u,\al)\|_p \leq \|u\|_p \|\al\|_p .
\end{equation}
Suppose that $(\al_\nu)\subset A\hot A$ is a bounded approximate diagonal for $A$ such that $\|\al_\nu\|_p \leq M$ for all $\nu$. We aim to show there is a bounded approximate diagonal $(U_\gamma) \subset C(X,A)\hot C(X,A)$ for $C(X,A)$ such that $\|U_\gamma\|_p \leq 2Mc$ for all $\gamma$, where $c>0$ is the constant asserted in Corollary~\ref{ineq}. To this end it suffices to show that, for any $\ep > 0$ and any finite set $F\subset C(X,A)$, there is $U= U_{(F, \ep)} \in C(X,A)\hot C(X,A)$ such that 
\[
\|U\|_p \leq 2Mc, \quad \|a\cdot U - U\cdot a\|_p < \ep \quad \text{and } \|\pi(U)a - a\|< \ep
\]
 for all $a\in F$. Indeed, for the natural partial order $(F_1,\ep_1) \prec (F_2,\ep_2)$ if and only if $F_1\subset F_2$ and $\ep_1\geq \ep_2$, the net $(U_{(F,\ep)})$ will be the desired approximate diagonal for $C(X,A)$.

Let $\ep >0$ and let $F\subset C(X,A)$ be a fixed finite set.

Case 1: assume each $a\in F$ is of the form $a = \sum_k {f_k a_k}$, where the sum is a finite sum, $f_k\in C(X)$ and $a_k\in A$. Let $N>0$ be an integer that is greater than the number of the terms of $a = \sum_k {f_k a_k}$ for all $a\in F$. Obviously, all elements $a_k$ associated to $a$ for all $a\in F$ form a finite set $F_A\subset A$, and all functions $f_k$ associated to $a$ for all $a\in F$ form a finite set $F_C \subset C(X)$. Let $L>0$ be a number such that $\|b\|_A \leq L$ for all $b\in F_A$ and $\|f\|_\infty \leq L$ for all $f\in F_C$. 

By the amenability of $A$, there is $\al \in (\al_\nu)$ such that
\[
\|b\cdot\al - \al\cdot b\|_p < \frac{\ep}{4cNL}, \quad \|\pi(\al)b-b\|_A < \frac{\ep}{NL} \quad (b\in F_A).
\]

On the other hand, by the compactness of $X$ there are finite open sets, say $V_i\subset X$ ($i=1,2,\ldots, n$), such that $X= \cup_i V_i$ and 
\[
|f(x) - f(y)| < \frac{\ep}{8c\|\al\|_p LN} \quad (f\in F_C, x,y\in V_i).
\]
From each $V_i$ we take a point $x_i$.
Apply partition of unity. We obtain continuous functions $h_1,h_2, \ldots, h_n \in C(X)$ such that Supp$(h_i)\subset V_i$, $0\leq h_i(x) \leq 1$ and $h_1+h_2+\cdots + h_n =1$ on $X$. Let $u_i = \sqrt{h_i}$ and set $u = \sum_{i=1}^n{u_i\otimes u_i}$. Then $u\in C(X)\otimes C(X)$ and $\pi(u) = 1$. From Corollary~\ref{ineq} it is not hard to see $\|u\|_p \leq 2c$ and
\[
\|f\cdot u - u\cdot f\|_p \leq \|\sum_i(f-f(x_i))u_i\otimes u_i\|_p + \|\sum_i u_i\otimes (f-f(x_i))u_i\|_p < \frac{\ep}{2\|\al\|_p LN}
\]
for all $f\in F_C$.

We now consider $U = T(u,\al)$. Then for $a = \sum_i {f_k a_k}\in F$ we have
\begin{align*}
\|a\cdot U- U\cdot a\|_p &= \|\sum_k{(T(f_ku,a_k\al) - T(uf_k,\al a_k))}\|_p\\
                   &= \|\sum_k{(T(f_ku,a_k\al-\al a_k) + T(f_k u - uf_k,\al a_k))}\|_p \\
                   &\leq \sum_k{(L\|u\|_p \|a_k\al-\al a_k\|_p + L\|\al\|_p \|f_k u - uf_k\|_p)} \\
                   &< NL\left(2c \frac{\ep}{4cNL} +\|\al\|_p \frac{\ep}{2\|\al\|_p LN}\right) =\ep;
\end{align*}
and
\begin{align*}
\|\pi(U)a -a\| &= \|\pi(u)\pi(\al)a -a\| =\|\sum_k{f_k\left(\pi(\al)a_k - a_k\right)}\|\\
               &\leq \sum_k{L\|\pi(\al)a_k - a_k\|_A} < NL \frac{\ep}{NL} =\ep.
\end{align*}
This completes the proof for case 1.

case 2: let $F$ be any finite set in $C(X,A)$.
We first observe that for $a\in C(X,A)$ and $\ep >0$ there is $a_\ep \in C(X,A)$ in the form $a_\ep = \sum_k {f_k a_k}$ such that $\|a - a_\ep\|_\infty < \ep$, where the right side of $a_\ep$ is a finite sum, $f_k\in C(X)$ and $a_k\in A$. In fact, since $X$ is compact, there are finite open sets $X_k\subset X$ such that $\cup_k X_k = X$ and $\|a(x)- a(y)\|_A < \ep$ for $x,y\in X_k$. Choose a point $x_k\in X_k$ for each $k$ and denote $a_k = a(x_k)$. From the partition of unity, there are $f_k \in C(X)$, such that supp$(f_k) \subset X_k$, $0\leq f_k(x) \leq 1$ for $x\in X$, and $\sum_k{f_k} = 1$. Then it is easy to check that $a_\ep = \sum_k{f_k a_k}$ satisfies the requirement.

From the above observation, for each $a\in F$ we may choose $a_\ep$ in the  form $a_\ep = \sum_k {f_k a_k}$ such that $\|a - a_\ep\|_\infty <\min\set{\ep/4, \ep/(8Mc)}$. Then $F_\ep = \set{a_\ep: a\in F}$ is a finite set of $C(X,A)$ satisfying the assumption of case 1. So there is $U\in C(X,A)\hot C(X,A)$ such that $\|U\|_p \leq 2cM$ and
\[
\|a_\ep\cdot U - U\cdot a_\ep\|_p < \ep/2, \quad \|\pi(U)a_\ep - a_\ep\| < \ep/2 \quad (a_\ep\in F_\ep).
\]
Then it can be checked easily that for this $U$
\[
\|a\cdot U - U\cdot a\|_p < \ep, \quad \|\pi(U)a-a\| < \ep \quad (a\in F).
\]
This completes the proof of the theorem.

\end{proof}

\begin{remark}
The converse of Theorem~\ref{bounded} is also true, i.e. if $C(X,A)$ is amenable then $A$ is amenable.
\end{remark}

\begin{proof}
Take an $x_0\in X$. Consider $T$: $C(X,A) \to A$ defined by $T(f) = f(x_0)$. This is a continuous subjective Banach algebra homomorphism. So $A$ is amenable if $C(X,A)$ is amenable.

\end{proof}

\section{Central compactly approximate diagonal}\label{unbounded a.d.}

In this section we consider when $C(X,A)$ has a compactly approximate diagonal. It is still unknown to the authors whether the existence of a compactly approximate diagonal for $A$ implies the existence of a compactly approximate diagonal for $C(X,A)$. We only have a partial answer to it.

\begin{prop}\label{cmp approx}
Let $X$ be a compact Hausdorff space and $A$ a Banach Algebra. If $A$ has a central compactly approximate diagonal, then $C(X,A)$ has a compactly approximate diagonal.
\end{prop}
\begin{proof}
It suffices to show that for every compact set $\K\subset C(X,A)$ and any $\ep>0$ there is $U\in C(X,A)\hot C(X,A)$ such that
\[
\|a\cdot U - U\cdot a\|_p < \ep, \quad \|\pi(U)a -a\|_\infty < \ep
\]
for all $a\in \K$.

Let $K=\set{a(x): a\in \K, x\in X}$. It is readily seen that $K$ is a compact set of $A$. Since $A$ has a central compactly approximate diagonal, there is $\al\in A\hot A$ such that 
\[
b\cdot \al - \al\cdot b = 0, \quad \|\pi(\al)b - b\|_A < \ep/2 \quad \text{for all } b\in K.
\]
 On the other hand, due to the compactness of $\K$ there are finite open sets $X_k\subset X$, $k=1,2,\ldots, n$, such that 
\[
\cup_k X_k = X, \quad \text{and } \|a(x)- a(y)\|_A < \frac{\ep}{8c\|\al\|_p} \quad (x,y\in X_k, a\in \K)
\]
 for each $k$, where $c>0$ is the number asserted in Corollary\ref{ineq}. Let $\set{f_1,f_2,\ldots,f_n} \subset C(X)$ be a partition of unity with respect to $\set{X_1,X_2,\ldots,X_n}$ (i.e. supp$(f_k) \subset X_k$, $0\leq f_k(x) \leq 1$ for $x\in X$, and $\sum_k{f_k} = 1$). Choose a point $x_k\in X_k$ for each $k$ and denote $a_k = a(x_k)$. Define, for each $a\in \K$,  $a_\ep = \sum_k{f_k a_k}$. Then
 \begin{equation}\label{approximate}
 \|a - a_\ep\|_\infty  = \|\sum_kf_k(a - a_k)\|_\infty \leq \frac{\ep}{8c\|\al\|_p} \quad (a\in \K).
 \end{equation}
 
 As we showed in the proof of \ref{bounded}, for the finite set $\set{f_1,f_2,\ldots,f_n}$ we have obtained, there is $u\in C(X)\otimes C(X)$ such that 
 \[
 \|u\|_p \leq 2c, \quad \pi(u) = 1, \quad \text{and } \|f_k\cdot u - u\cdot f_k\|_p < \frac{\ep}{2nM\|\al\|_p} \quad (k=1,2,\ldots, n)
 \]
 where $M= \sup\set{\|b\|_A : b\in K}$.
 
 We let $U = T(u,\al)$. Then $\|U\|_p \leq 2c\|\al\|_p$. For each $a\in \K$ we first have
 \begin{align*}
\|a_\ep\cdot U- U\cdot a_\ep\|_p &= \|\sum_k{(T(f_ku,a_k\al) - T(uf_k,\al a_k))}\|_p\\
                   &= \|\sum_k{(T(f_ku,a_k\al-\al a_k) + T(f_k u - uf_k,\al a_k))}\|_p \\
                   &= \|\sum_k{T(f_k u - uf_k,\al a_k)}\|_p\\
                   &\leq \sum_k{\|f_k u - uf_k\|_p\|a_k\|_A \|\al\|_p}  \\
                   &< \frac{\ep}{2nM\|\al\|_p}nM\|\al\|_p  =\ep/2.
\end{align*}
and
\begin{align*}
\|\pi(U)a_\ep -a_\ep\| &= \|\pi(u)\pi(\al)a_\ep -a_\ep\|\\ 
                       &=\|\sum_k{f_k\left(\pi(\al)a_k - a_k\right)}\|
               \leq \frac{\ep}{2}\norm{\sum_k{f_k}}_\infty  =\ep/2.
\end{align*}

Without loss of general we may assume $\|\al\|_p \geq 1$. Then by inequality~(\ref{approximate}) the above estimates immediately lead to
\[
\|a\cdot U - U\cdot a\|_p < \ep \quad \text{and } \|\pi(U)a -a\|_\infty < \ep \quad (a\in \K).
\]
The proof is complete.

\end{proof}

\section{Weak amenability}\label{weak amen}
To conclude this note we address a brief discussion in this section to the weak amenability of $C(X,A)$.
We recall that a Banach algebra $A$ is weakly amenable if every continuous derivation from $A$ into $A^*$ is inner. It is unknown whether $C(X,A)$ is weakly amenable if $A$ is so, except for the case that $A$ is a C*-algebra. In the latter case $C(X,A)$ itself is a C*-algebra. We remark that in general, the dual space $C(X,A)^*$ of $C(X,A)$ is identical to $I_1\left(C(X),A^*\right)$, the space of all 1-integral operators from $C(X)$ into $A^*$ \cite[page 18]{Pisier}. 

For the commutative case we have a positive answer to the above question. According to \cite{B-C-D} a commutative Banach algebra $A$ is weakly amenable if and only if every continuous derivation from $A$ into a commutative Banach $A$-bimodule $E$ is trivial. Here a Banach $A$-bimodule $E$ is a commutative bimodule if $a\cdot x =x\cdot a$ for all $a\in A$ and $x\in E$.

\begin{prop}
Let $X$ be a compact Hausdorff space and $A$ a commutative Banach algebra. If $A$ is weakly amenable and if $A$ has a bounded approximate identity, then $C(X,A)$ is weakly amenable.
\end{prop}

\begin{proof}
Clearly, $C(X,A)$ is a commutative Banach algebra and $A$ is a closed subalgebra of it. Therefore, $C(X,A)$ is naturally a commutative $A$-bimodule, and it is also a commutative $C(X)$-bimodule. Let $D$: $C(X,A) \to C(X,A)^*$ be a continuous derivation. Then $D|_A$: $A\to C(X,A)^*$ is a continuous derivation. Since $A$ is commutative and weakly amenable, $D|_A = 0$.

Let $(e_\nu)$ be a bounded approximate identity of $A$. Then $(e_\nu)$ is also a bounded approximate identity for $C(X,A)$ and {wk*}-$\lim D(e_\nu) = 0$. For each $f\in C(X)$, we note that wk*-$\lim D(fe_\nu)$ exists. To see this it suffices to show that all weak* convergent subnets of $(D(fe_\nu))$ converge to the same limit. Assume wk*-$\lim D(fe_i)$ and wk*-$\lim D(fe_j)$ exist, where $(e_i), (e_j)$ are subnets of $(e_\nu)$. Then
\begin{align*}
D(fe_i) = \lim_j D(fe_je_i) &= \lim_j (fe_j)D(e_i) + D(fe_j)e_i \\
                            &=f D(e_i) + \text{wk*-}\lim_j D(fe_j)e_i .
\end{align*}
Take weak* limit in $i$. We then get 
\[
\text{wk*-}\lim D(fe_i)= \text{wk*-}\lim D(fe_j)
\]
 as claimed. So $\tilde{D}$: $C(X)\to C(X,A)^*$ given by $\tilde{D}(f) = \text{wk*-}\lim D(fe_\nu)$ is well-defined. 
\begin{align*}
\tilde{D}(fg) &=\text{wk*-}\lim D(fge_\nu) = \text{wk*-}\lim_\nu (\lim_\mu D(fe_\mu ge_\nu))\\
              & = \text{wk*-}\lim_\nu fD( ge_\nu) + \text{wk*-}\lim_\mu D(fe_\mu) g =
f\tilde{D}(g) + \tilde{D}(f)g
\end{align*}
for all $f,g\in C(X)$.
Therefore $\tilde{D}$ is a derivation. However, $C(X)$ is amenable. We have $\tilde{D} = 0$. Now for $f\in C(X)$ and $a\in A$, we have
\[
D(fa) = \tilde{D}(f)\cdot a + f D_A(a) =0.
\]
So $D =0$ on the linear span of $\set{fa: f\in C(X), a\in A}$. On the other hand every element of $C(X,A)$ may be approached in norm by elements from this linear span as we have seen in the proofs of Theorem~\ref{bounded} and Proposition~\ref{cmp approx}. We thus derive $D=0$ on the whole $C(X,A)$.
The proof is complete.
\end{proof}

\end{document}